\documentclass[12pt,reqno]{amsart}

\newcommand\version{July 21, 2024}

\usepackage{amsmath,amsfonts,amsthm,amssymb,amsxtra}
\usepackage{bbm} 
\usepackage{hyperref} 

\newtheorem{theorem}{Theorem}

\theoremstyle{definition}

\theoremstyle{remark}

\newcommand{\const}{\mathrm{const}\, }

\renewcommand{\epsilon}{\varepsilon}

\renewcommand{\phi}{\varphi}
\newcommand{\R}{\mathbb{R}}

\newcommand{\Z}{\mathbb{Z}}

\DeclareMathOperator{\Tr}{Tr}

\begin{document}

\title[Minimizing Schr\"odinger eigenvalues for confining potentials]{Minimizing Schr\"odinger eigenvalues\\ for confining potentials}

\author{Rupert L. Frank}
\address[Rupert L. Frank]{Mathe\-matisches Institut, Ludwig-Maximilians Universit\"at M\"un\-chen, The\-resienstr.~39, 80333 M\"unchen, Germany, and Munich Center for Quantum Science and Technology, Schel\-ling\-str.~4, 80799 M\"unchen, Germany, and Mathematics 253-37, Caltech, Pasa\-de\-na, CA 91125, USA}
\email{r.frank@lmu.de}

\thanks{\emph{Date:} \version \\
	\copyright\, 2024 by the author. This paper may be reproduced, in its entirety, for non-commercial purposes.\\
Partial support through US National Science Foundation (DMS-1954995), as well as through the German Research Foundation (EXC-2111-390814868 and TRR 352-470903074) is acknowledged.}

\begin{abstract}
	We consider the problem of minimizing the lowest eigenvalue of the Schr\"odinger operator $-\Delta+V$ in $L^2(\R^d)$ when the integral $\int e^{-tV}\,dx$ is given for some $t>0$. We show that the eigenvalue is minimal for the harmonic oscillator and derive a quantitative version of the corresponding inequality.
\end{abstract}

\maketitle

\subsection*{Introduction and main result}

We are interested in a sharp lower bound on the lowest eigenvalue $E_0$ of the Schr\"odinger operators $-\Delta+V$ in $L^2(\R^d)$, $d\geq 1$, with a real confining potential $V$, and in the stability properties of the resulting functional inequality.

More precisely, we assume that $e^{-tV}\in L^1(\R^d)$ for some $t>0$. Under this assumption the operator $-\Delta+V$ has compact resolvent as a consequence of the Golden--Thompson inequality; see \cite[Theorem 9.2]{Si} and \eqref{eq:gt} below.

The inequality that we are interested in asserts that the lowest eigenvalue $E_0$ satisfies
\begin{equation}
	\label{eq:evineq}
	\boxed{E_0 \geq - t^{-1} \ln\int_{\R^d} e^{-t V(x)}\,dx + dt^{-1} (1+ \tfrac12\ln(\pi t)) \,.}
\end{equation}
While we have not found this inequality stated explicitly in the literature, it is probably known to experts in the field. We will review its proof after the statement of our main theorem.

Inequality \eqref{eq:evineq} is sharp for any given $t>0$. Indeed, equality holds if
\begin{equation}
	\label{eq:equality}
	V(x) = t^{-2} |x-b|^2 + \const
\end{equation}
for some $b\in\R^d$, as follows from the explicit knowledge of the eigenvalues of the harmonic oscillator; see, e.g., \cite[Subsection 4.2.1]{FrLaWe}. Conversely, one can show that if equality in \eqref{eq:evineq} holds, then $V$ is necessarily of the form \eqref{eq:equality}. We will show this below after having given the proof of \eqref{eq:evineq}.

Our main result is a stability version of this characterization of optimizers. Namely, we shall show that if $V$ almost saturates \eqref{eq:evineq}, then $V$ is almost of the form \eqref{eq:equality}. We prove this approximate form of $V$ in the natural norm of the problem and with the optimal exponent. The precise statement is as follows.

\begin{theorem}\label{main}
	There is a constant $c>0$ such that for all $d\geq 1$, all $t>0$ and all $V$ with $e^{-tV}\in L^1(\R^d)$,
	\begin{equation*}
		\boxed{
		\begin{aligned}
			E_0 & \geq - t^{-1} \ln\int_{\R^d} e^{-t V(x)}\,dx + dt^{-1} (1+ \tfrac12\ln(\pi t)) \\
			& \quad + c \, t^{-1} \inf_{b\in\R^d} \left\|  \frac{e^{-t V}}{\int_{\R^d} e^{-t V}\,dx} - (\pi t)^{-\frac d2} e^{-\frac{1}{t}|x-b|^2} \right\|_{L^1(\R^d)}^2.
		\end{aligned}
		}
	\end{equation*}
\end{theorem}

We emphasize that the constant $c$ is independent of $d$, $t$ and $V$. In fact, this constant is explicit in the sense that no compactness argument is used in its proof.


\subsection*{Background}

Before turning to the proofs of \eqref{eq:evineq} and of Theorem \ref{main}, let us put this inequality and its stability result into perspective. The problem of minimizing $E_0$ among potentials $V$ with fixed $L^p$ norm was proposed by Keller \cite{Ke} and solved by Lieb and Thirring \cite{LiTh} in terms of optimal constants in certain Sobolev interpolation inequalities. In this way one is led to the functional inequality
\begin{equation}
	\label{eq:lt1}
	E_0 \geq - \left( L_{\gamma,d}^{(1)} \int_{\R^d} V(x)_-^{\gamma+\frac d2}\,dx \right)^\frac 1\gamma
\end{equation}
for a certain constant $L_{\gamma,d}^{(1)}$, provided $\gamma\geq\frac12$ if $d=1$ and $\gamma>0$ if $d\geq 2$; see \cite[Subsection 5.1.2]{FrLaWe}. Here $v_-:=\max\{-v,0\}$ denotes the negative part.

Similarly, our inequality \eqref{eq:evineq} can be considered as the answer to the problem of minimizing $E_0$ among potentials $V$ for which the integral $\int e^{-tV}\,dx$ is fixed for some $t>0$. 

The stability of the functional inequality \eqref{eq:lt1} was investigated in \cite{CaFrLi}, based on two stability results for a Sobolev interpolation inequality and for H\"older's inequality. Similarly, our proof of Theorem \ref{main} is based on a stability result for the log-Sobolev inequality and for a certain endpoint case of H\"older's inequality. The relevant stability result for the log-Sobolev inequality was recently shown in \cite{DoEsFiFrLo} (see also \cite{DoEsFiFrLo2} for a partially alternative proof) and that for the endpoint case of H\"older's inequality is the Csisz\'ar--Kullback--Pinsker inequality in the form of Carlen in~\cite{Ca2}. Theorem \ref{main} will follow by combining these two results. We emphasize that there is a large literature on stability versions of the log-Sobolev inequality (see the references in \cite{DoEsFiFrLo}), but it is the new form of the remainder term in \cite{DoEsFiFrLo} that is required for our application here. Also the use of the Csisz\'ar--Kullback--Pinsker inequality for the stability of functional inequalities involving exponential integrability is very recent. It was used by Carlen in~\cite{Ca2} to quantify the rate of approach to equilibrium for the critical mass Keller--Segel system, improving earlier results in \cite{CaFi}. This is based on a stability result for the logarithmic Hardy--Littlewood--Sobolev inequality.

The question of whether there is a quantitative version of the sharp Sobolev inequality was raised by Brezis and Lieb in \cite{BrLi}, where also initial results were obtained. It was answered by Bianchi and Egnell in \cite{BiEg} using a technique that has proved to be very influential. In the last decade there has been an enormous interest in variants of this question and we refer to \cite{ChFrWe,CaFi,Ca1,DoEsFiFrLo,Fr,Ca2} for further references. Concerning eigenvalue inequalities for Schr\"odinger operators we also mention the stability results in \cite{Ma,MaRB}, which were obtained by techniques different from ours.

As an aside, we mention that, just like \eqref{eq:lt1} is the `one-particle' version of the Lieb--Thirring inequality (see \cite[Subsection 5.1.2]{FrLaWe} and \cite{Fr0}), inequality \eqref{eq:evineq} is the `one-particle' version of the Golden--Thompson inequality
\begin{equation}
	\label{eq:gt}
	\Tr e^{t(\Delta-V)} \leq (4\pi t)^{-\frac d2} \int_{\R^d} e^{-tV(x)}\,dx \,;
\end{equation}
see, e.g., \cite[Theorem 9.2]{Si}. Indeed, if we bound the left side from below by $e^{-tE_0}$, we get an inequality of the form \eqref{eq:evineq}, but with the term $dt^{-1} (1+ \tfrac12\ln(\pi t))$ replaced by the smaller term $\tfrac d2 t^{-1} \ln(4\pi t)$. Thus, this does not yield a sharp inequality. The Golden--Thompson inequality is sharp in the semiclassical regime of many eigenvalues. A similar phenomenon occurs for the entropy inequality in \cite{FrLi}, obtained from a more general form of the Golden--Thompson inequality, which is the `many-particle' version of the Beckner--Hirschman inequality \cite{Hi,Be}. It would be interesting to obtain a quantitative version of the latter inequality.


\subsection*{Acknowledgement}
The author is very grateful to Eric Carlen for several useful conversations on the topic of this paper and on \cite{Ca1,Ca2}.


\subsection*{Proof of \eqref{eq:evineq}}

As a warm-up, we show how the lower bound on $E_0$ in \eqref{eq:evineq} can be derived from the log-Sobolev inequality and the Gibbs variational principle. Each step in this proof will later be made quantiative in the proof of Theorem \ref{main}.

We shall use the variational characterization of the lowest eigenvalue (see, e.g., \cite[Corollary 1.11]{FrLaWe}), viz. 
\begin{align}\label{eq:varprinc}
	E_0 & = \inf_{\|\psi\|=1} \int_{\R^d} \left( |\nabla\psi|^2 +V\psi^2\right)dx \,.
\end{align}
Here $\|\cdot\|$ denotes the $L^2$-norm. We bound the gradient term from below by means of the (Euclidean) log-Sobolev inequality \cite[Theorem 8.14]{LiLo}, which says that
\begin{equation}
	\label{eq:logsob}
	\lambda^2 \int_{\R^d} |\nabla\psi|^2\,dx  - \pi \int_{\R^d} |\psi|^2\ln(|\psi|^2)\,dx \geq d\pi (1+\ln\lambda)
\end{equation}
for any $\psi\in H^1(\R^d)$ with $\|\psi\| =1$ and any $\lambda>0$. Consequently,
\begin{align}\label{eq:evlowerlogsob}
	E_0 \geq \inf_{\|\psi\|=1} \int_{\R^d} \left(V + \pi\lambda^{-2}\ln(|\psi|^2)\right) |\psi|^2 \,dx + d\pi \lambda^{-2}(1+\ln\lambda) \,.
\end{align}
A simple and well-known computation, which is known as Gibbs variational principle in statistical mechanics, shows that on any measure space $X$ (with measure denoted by $dx$)
\begin{equation}
	\label{eq:gibbs}
	\inf_{\|\rho\|_1=1,\, \rho\geq 0} \int_X ( V+ T \ln \rho)\rho\,dx = -T \ln\int_X e^{-\frac1T V}\,dx \,.
\end{equation}
Here $\|\cdot\|_1$ denotes the $L^1$-norm. Inserting the lower bound given by \eqref{eq:gibbs} with $T=\pi\lambda^{-2}=t^{-1}$ and $\rho = |\psi|^2$ into \eqref{eq:evlowerlogsob}, we obtain \eqref{eq:evineq}.

In the same way, we can characterize the cases of equality in \eqref{eq:evineq}. As shown by Carlen \cite{Ca}, equality in \eqref{eq:logsob} holds only if $\psi(x) = e^{i\theta} \lambda^{-\frac d2} e^{-\frac\pi{2\lambda^2} |x-b|^2}$ for some $b\in\R^d$, $\lambda>0$ and $\theta\in\R/2\pi\Z$. Moreover, the infimum in \eqref{eq:gibbs} is attained if and only if $\rho = Z^{-1} e^{-\frac1TV}$ with $Z=\int_X e^{-\frac1T V}\,dx$. Thus, equality in \eqref{eq:evineq} forces $V(x) = t^{-2} |x-b|^2 - t^{-1} \ln(Z/(\pi t)^{d/2})$, as claimed.

In the above argument we were somewhat cavalier about finiteness of certain integrals. While for the proof of inequality \eqref{eq:evineq} we can restrict ourselves to functions $\psi\in C^1_c(\R^d)$, say, for which all manipulations are justified, for the characterization of cases of equality we need to work with general function $\psi$ in the form domain of the Schr\"odinger operator $-\Delta+V$. We defer the discussion of the relevant technical details to the end of this paper.


\subsection*{Proof of Theorem \ref{main}}

	The stability version of the log-Sobolev inequality \eqref{eq:logsob} from \cite{DoEsFiFrLo} reads
	\begin{align*}
		& \lambda^2 \int_{\R^d} |\nabla\psi|^2\,dx  - \pi \int_{\R^d} |\psi|^2\ln(|\psi|^2)\,dx - d\pi (1+\ln\lambda) \\
		& \geq \kappa \inf_{b\in\R^d,\, \theta\in\R/2\pi\Z} \int_{\R^d} |\psi - e^{i\theta} \lambda^{-\frac d2} e^{-\frac{\pi}{2\lambda^2}|x-b|^2}|^2\,dx
	\end{align*}
	for all $\psi\in H^1(\R^d)$ with $\|\psi\|=1$ with a universal constant $\kappa>0$ (independent of $d$). (The inequality there is stated for Gaussian measure, but it can easily be transformed into the above form with Lebesgue measure. The dependence on $\lambda$ follows by scaling.) We shall apply this inequality only to nonnegative functions $\psi$, in which case the infimum over $\theta$ is attained when $e^{i\theta}=1$, so the inequality becomes
	\begin{align}\label{eq:logsobstab}
		& \lambda^2 \int_{\R^d} |\nabla\psi|^2\,dx  - \pi \int_{\R^d} \psi^2\ln(\psi^2)\,dx - d\pi (1+\ln\lambda) \notag \\
		& \geq \kappa \inf_{b\in\R^d} \int_{\R^d} (\psi - \lambda^{-\frac d2} e^{-\frac{\pi}{2\lambda^2}|x-b|^2})^2\,dx \,.
	\end{align}
	
	Meanwhile, the Csisz\'ar--Kullback--Pinsker inequality in Carlen's reformulation as a stability version of the Gibbs variational principle \eqref{eq:gibbs} reads
	\begin{equation}
		\label{eq:gibbsstab}
			\int_X ( V+ T \ln \rho)\rho\,dx + T \ln\int_X e^{-\frac1T V}\,dx \geq \frac T2 \left\| \rho - \frac{e^{-\frac1T V}}{\int_X e^{-\frac1T V}\,dx} \right\|_1^2
	\end{equation}
	for all functions $\rho\geq 0$ on $X$ with $\|\rho\|_1=1$; see \cite[Theorem 3.2]{Ca2} and our discussion at the end of this paper.
	
	Let $\psi_0$ is a normalized eigenfunction corresponding to $E_0$ and note that
	\begin{equation}
		\label{eq:gsenergy}
		E_0 = \int_{\R^d} \left( |\nabla\psi|^2 +V|\psi|^2\right)dx \,.
	\end{equation}
	It is well known that $\psi_0$ may be chosen nonnegative. Combining inequalities \eqref{eq:logsobstab} and \eqref{eq:gibbsstab} with \eqref{eq:gsenergy} and choosing as before $T=\pi\lambda^{-2}=t^{-1}$ gives
	\begin{align*}
		& E_0 + t^{-1} \ln\int_{\R^d} e^{-t V(x)}\,dx - dt^{-1} (1+ \tfrac12\ln(\pi t)) \\
		& \geq (\pi t)^{-1} \kappa \inf_{b\in\R^d} \int_{\R^d} (\psi_0 - (\pi t)^{-\frac d4} e^{-\frac{1}{2t}|x-b|^2})^2\,dx + (2t)^{-1} \left\| \psi_0^2 - \frac{e^{-t V}}{\int_{\R^d} e^{-t V}\,dx} \right\|_1^2.
	\end{align*} 
	
	We recall that for measurable functions $f,g\geq 0$ on a measure space $X$ with $\int_X f\,dx = \int_X g\,dx =1$ one has
	$$
	\|f- g\|_1 \leq 2\, \|\sqrt f-\sqrt g\|_2 \,.
	$$
	Indeed, on the left side we write $f-g=(\sqrt f-\sqrt g)(\sqrt f+\sqrt g)$ and we apply Schwarz's inequality. The assertion follows since $\|\sqrt f+\sqrt g\|_2 \leq \|\sqrt f\|_2 + \|\sqrt g\|_2 = 2$.
	
	In our situation this implies that
	$$
	\int_{\R^d} (\psi_0 - (\pi t)^{-\frac d4} e^{-\frac{1}{2t}|x-b|^2})^2\,dx
	\geq \frac14 \left( \int_{\R^d} \left|\psi_0^2 - (\pi t)^{-\frac d2} e^{-\frac{1}{t}|x-b|^2}\right| dx \right)^2.
	$$
	Thus, we find
	\begin{align*}
		& E_0 + t^{-1} \ln\int_{\R^d} e^{-t V(x)}\,dx - dt^{-1} (1+ \tfrac12\ln(\pi t)) \\
		& \geq (4\pi t)^{-1} \kappa \inf_{b\in\R^d} \left\| \psi_0^2 - (\pi t)^{-\frac d2} e^{-\frac{1}{t}|x-b|^2} \right\|_1^2 + (2t)^{-1} \left\| \psi_0^2 - \frac{e^{-t V}}{\int_{\R^d} e^{-t V}\,dx} \right\|_1^2 \\
		& \geq (2t)^{-1} \min\{\tfrac\kappa{2\pi},1\} \inf_{b\in\R^d} \left( \left\| \psi_0^2 - (\pi t)^{-\frac d2} e^{-\frac{1}{t}|x-b|^2} \right\|_1^2 + \left\| \psi_0^2 - \frac{e^{-t V}}{\int_{\R^d} e^{-t V}\,dx} \right\|_1^2 \right) \\
		& \geq (4t)^{-1} \min\{\tfrac\kappa{2\pi},1\} \inf_{b\in\R^d} \left( \left\| \psi_0^2 - (\pi t)^{-\frac d2} e^{-\frac{1}{t}|x-b|^2} \right\|_1 + \left\| \psi_0^2 - \frac{e^{-t V}}{\int_{\R^d} e^{-t V}\,dx} \right\|_1 \right)^2 \\
		& \geq (4t)^{-1} \min\{\tfrac\kappa{2\pi},1\} \inf_{b\in\R^d} \left\|  \frac{e^{-t V}}{\int_{\R^d} e^{-t V}\,dx} - (\pi t)^{-\frac d2} e^{-\frac{1}{t}|x-b|^2} \right\|_1^2.
	\end{align*}
	This proves the claimed inequality with $c=4^{-1} \min\{\tfrac\kappa{2\pi},1\}$.	
\qed


\subsection*{Some technicalities}

Let us discuss some technical details when carrying out the proof of \eqref{eq:equality} and, in particular, when characterizing the cases of equality. For the latter purpose, we need to carry out the proof for arbitrary functions $\psi$ in the form domain of the Schr\"odinger operator $-\Delta+V$. In the proof the integrals $\int V|\psi|^2\,dx$ and $\int |\psi|^2\ln|\psi|^2\,dx$ appear and our goal now is to show that both integrals converge absolutely.

The elementary inequality $w^p \leq (p/e)^p e^w$ for every $w\geq 0$ and $p>0$ implies that $\int_{\R^d} V_-^p\,dx \leq (p/te)^p \int_{\{V<0\}} e^{-tV}\,dx$. Applying this with $p>\max\{\frac d2,1\}$, we see that $V_-$ is infinitesimally form bounded with respect to $-\Delta$; see \cite[Proposition 4.3]{FrLaWe}. In particular, for functions in the form domain all three integrals $\int |\nabla\psi|^2\,dx$, $\int V_+|\psi|^2\,dx$ and $\int V_-|\psi|^2\,dx$ are finite. Using Sobolev inequalies, we easily see that $\int |\psi|^2 \ln_+ |\psi|^2\,dx <\infty$. Note that at this point we do not know yet that $\int |\psi|^2 \ln_- |\psi|^2\,dx$ is finite.

We use the Gibbs variational principle in the form of the nonnegativity of the relative entropy
\begin{equation}
	\label{eq:relent}
	\int_X \rho (\ln\rho - \ln\rho_0)\,dx \geq 0 \,,
\end{equation}
valid for all measurable functions $\rho,\rho_0\geq 0$ on $X$ with $\int_X \rho\,dx=\int_X \rho_0\,dx<\infty$. By definition, the left side of \eqref{eq:relent} is $+\infty$ if $\int_{\{\rho_0=0\}} \rho\,dx >0$. Note that the integral \eqref{eq:relent} is well defined, though possibly $+\infty$, since in view of the inequality $\alpha \ln(1/\alpha)\leq 1/e$ for $\alpha\in(0,1]$ we have
$$
\int_X \rho (\ln\rho - \ln \rho_0)_- \,dx = \int_{\{\rho<\rho_0\}} \rho_0 \frac{\rho}{\rho_0} \ln\frac{\rho_0}{\rho}\,dx \leq e^{-1} \int_{\{\rho<\rho_0\}} \rho_0 <\infty \,.
$$

For the sake of completeness let us give the proof of \eqref{eq:relent}. Concavity of the function $p\mapsto - p \ln p=:s(p)$ implies that $s(p)-s(q)\leq s'(q)(p-q)$, so pointwise on $X$ we have
$$
-\rho(x)\ln\rho(x) + \rho_0(x)\ln\rho_0(x) \leq - (1+\ln\rho_0(x))(\rho(x)-\rho_0(x)) \,,
$$
that is,
$$
\rho(x) \ln\rho(x) - \rho(x)\ln\rho_0(x) \geq \rho(x)-\rho_0(x) \,.
$$
Integrating this inequality, recalling integrability of the negative part of the left side, gives \eqref{eq:relent}.

We apply inequality \eqref{eq:relent} with $X=\R^d$, endowed with Lebesgue measure, and with $\rho=|\psi|^2$ and $\rho_0 = Z^{-1} e^{-tV}$, where $Z=\int e^{-tV}\,dx$. Here $\psi$ belongs to the form domain of $-\Delta+V$ and is normalized by $\|\psi\|=1$. According to the properties of functions in the form domain, the integral $\int \rho \ln\rho_0\,dx$ converges absolutely and $\int \rho\ln_+\rho\,dx<\infty$. In view of \eqref{eq:relent}, this implies that $\int \rho\ln_-\rho\,dx<\infty$. Therefore all integrals appearing in the proof of \eqref{eq:evineq} converge absolutely and all manipulations are justified.


\subsection*{On inequality \eqref{eq:gibbsstab}}

An inequality of Pinsker, in the sharp form of Csisz\'ar \cite{Cs} and Kullback \cite{Ku}, states that
\begin{equation}
	\label{eq:ckp}
	\int_X \rho (\ln\rho - \ln\rho_0)\,dx \geq \frac12 \left\| \rho - \rho_0 \right\|_1^2 \,,	
\end{equation}
valid for all measurable functions $\rho,\rho_0\geq 0$ on $X$ with $\int_X \rho\,dx=\int_X \rho_0\,dx=1$. As before, we use the convention that the left side of \eqref{eq:ckp} is $+\infty$ if $\int_{\{\rho_0=0\}} \rho\,dx >0$. Note that \eqref{eq:ckp} is a quantitative version of the nonnegativity of the relative entropy in \eqref{eq:relent}.

Choosing, similarly as before, $X=\R^d$, endowed with Lebesgue measure, $\rho=|\psi|^2$ and $\rho_0 = Z^{-1} e^{-\frac1TV}$, where $Z=\int e^{-\frac1TV}\,dx$, we obtain \eqref{eq:gibbsstab}. The necessary manipulations are all legitimate due to the absolute convergence of the integrals that we have just discussed.

Several proofs of \eqref{eq:ckp} are known; see, e.g., \cite[Proof of Theorem 22.10, Exercise 29.22, as well as the bibliographic notes to Chapter 22]{Vi}. A recent one was given by Carlen in \cite[Theorem 3.2]{Ca2}, based on \cite[Theorem 2.12]{Ca1}. Inspection of the proof shows that the assumption there that $dx$ is a probability measure is not needed.

\bibliographystyle{amsalpha}

\end{document}